\documentclass[english, 11pt]{amsart}

\usepackage{tikz}
\usepackage{aliascnt}
\usepackage{mathrsfs}
\usepackage[all,poly,knot]{xy}
\usepackage{hyperref}
\usepackage{comment}  
\usepackage{csquotes}   
\usepackage{amssymb,amsbsy,amsmath,amsfonts,amssymb,amscd,
	graphics,color,footmisc,fancyhdr,multicol,fancybox,
	graphicx,mathrsfs,rotating,ifthen,wasysym}

\usepackage{times}
\usepackage{pdfpages}
\usepackage{multirow}
\usepackage{nccrules}
\usepackage{textcomp}
\usepackage{comment}
\usepackage{framed}
\usepackage[all]{xy}
\usepackage{graphicx}
\usepackage{pgf,tikz}
\usepackage{mathrsfs}
\usetikzlibrary{arrows}
\usepackage{lipsum}
\usepackage{mathtools}

\DeclareMathOperator{\dif}{\text{\normalfont d}}

\DeclareMathOperator{\rank}{rank}

\DeclareMathOperator{\supp}{supp}
\DeclareMathOperator{\ord}{ord}

\newcommand{\BB}{\mathbb{B}}
\def\log{\mathrm{log}\,}

\theoremstyle{plain}

\newtheorem{thm}{Theorem}[section]  

\newtheorem{cor}[thm]{{Corollary}}

\newtheorem{rem}[thm]{Remark}
\newtheorem{defi}[thm]{Definition}
\newtheorem*{claim}{Claim}

\theoremstyle{remark}

\numberwithin{equation}{section}

\usepackage[top=1.5in, bottom=1.5in, left=1in, right=1in]{geometry}

\usepackage[hyperpageref]{backref}

\usepackage{xcolor}
\hypersetup{
	colorlinks,
	linkcolor={red!50!black},
	citecolor={blue!62!black},
	urlcolor={blue!80!black}
}
\theoremstyle{plain}
\newcommand{\thistheoremname}{}
\newtheorem*{genericthm*}{\thistheoremname}
\newenvironment{namedthm*}[1]{\renewcommand{\thistheoremname}{#1}%
	\begin{genericthm*}}
	{\end{genericthm*}}

\newtheoremstyle{named}{}{}{\itshape}{}{\bfseries}{.}{.5em}{\thmnote{#3's }#1}
\theoremstyle{named}

\makeatletter
\newcommand\thankssymb[1]{\textsuperscript{\@fnsymbol{#1}}}
\makeatother
\begin{document} 
	\title[Degeneracy of holomorphic mappings into or avoiding Fermat type hypersurfaces]{\bf Degeneracy of holomorphic mappings\\
		into or avoiding Fermat type hypersurfaces}

	
	\author{Dinh Tuan Huynh}
	
	\address{Department of Mathematics, University of Education, Hue University, 34 Le Loi St., Hue City, Vietnam}
	\email{huynhdinhtuan@dhsphue.edu.vn}

\begin{abstract}
We prove that if $f\colon\mathbb{C}^p\rightarrow\mathbb{P}^n(\mathbb{C})$ is a holomorphic mapping of maximal rank whose image lies in the Fermat hypersurface of degree $d>(n+1)\max\{n-p,1\}$, then its image is contained in a linear subspace of dimension at most $\bigg[\dfrac{n-1}{2}\bigg]$. Analog in the logarithmic case is also given. Our result strengthens a classical result of Green and provides a Nevanlinna theoretic proof for a recent result due to Etesse.
\end{abstract}
\maketitle

\section{Introduction}

In complex analysis, the Little Picard Theorem states that any meromorphic maps on the complex plane avoiding $3$ distinct points must be constant. Higher dimensional generalizations of this classical result expect the degeneracy of holomorphic mappings from $\mathbb{C}^p$ into a projective variety $X$ omitting a  hypersurface $D$ of large enough degree. Green and Griffiths \cite{Green75,Green-Griffiths1980} anticipated that in the case where $X=\mathbb{C}\mathbb{P}^n$, the expected lower degree bound of $D$ should be $n+2$.
\begin{namedthm*}{Green-Griffiths'conjecture (logarithmic case)}
If $D$ is a simple normal crossing divisor on the projective space $\mathbb{C}\mathbb{P}^n$ of degree $d\geq n+2$, then the image of any holomorphic mapping $f\colon\mathbb{C}^p\rightarrow\mathbb{C}\mathbb{P}^n$ omitting $D$ lies in some proper algebraic subvariety of $\mathbb{C}\mathbb{P}^n$.	
\end{namedthm*}

In the so-called {\it compact} case, they also predicted the following

\begin{namedthm*}{Green-Griffiths'conjecture (compact case)}
	If $D$ is a nonsingular hypersurface in the projective space $\mathbb{C}\mathbb{P}^n$ of degree $d\geq n+2$, then the image of any holomorphic mapping $f\colon\mathbb{C}^p\rightarrow D$ lies in some proper algebraic subvariety of $D$.	
\end{namedthm*}

Many works have been done during recent decades towards the above conjecture. Notably, Green-Griffiths' conjecture was confirmed under the condition that the degree of $D$ is very high compared with the dimension. Currently the effective lower degree bound is still very far from the expected optimal one (exponential growth compared with the dimension). The readers are referred to \cite{Siu2004,DMR2010,Demailly2012,Demailly2020,Berczi2019,Darondeau2016,HVX2019} and references there in for recent progresses around this problem.

Green himself tested the above conjecture in various situations \cite{Green75}. In the case where $D$ is a Fermat type hypersurface, he obtained the following results.
\begin{thm}[compact case]
	Let $p\geq 1,n\geq2$ be positive integers. Let $F$ be the Fermat hypersurface of degree $d$ in $\mathbb{C}\mathbb{P}^n$, defined by the homogeneous polynomial 
$$
Q(\omega)=\sum_{i=0}^n\omega_i^d,
$$
where $\omega=[\omega_0:\omega_1:\dots:\omega_n]$ is a homogeneous coordinate of $\mathbb{C}\mathbb{P}^n$.
If $d>n^2-1$, then the image of every holomorphic map $f\colon\mathbb{C}^p\rightarrow F$ lies in a linear subspace of dimension at most $\bigg[\dfrac{n-1}{2}\bigg]$.	
\end{thm}

\begin{thm}[Logarithmic case]
Let $p\geq 1,n\geq2$ be positive integers. Let $F$ be the Fermat hypersurface of degree $d$ in $\mathbb{C}\mathbb{P}^n$, defined by the homogeneous polynomial 
$$
Q(\omega)=\sum_{i=0}^n\omega_i^d,
$$
where $\omega=[\omega_0:\omega_1:\dots:\omega_n]$ is a homogeneous coordinate of $\mathbb{C}\mathbb{P}^n$.
If $d>n(n+1)$, then the image of every holomorphic map $f\colon\mathbb{C}^p\rightarrow\mathbb{C}\mathbb{P}^n$  omitting the Fermat hypersurface $F$ lies in a linear subspace of dimension at most $\bigg[\dfrac{n}{2}\bigg]$.	
\end{thm}

In this paper, we improve the above two theorems for the class of holomorphic map of maximal rank,  by decreasing the lower degree bound.

\begin{defi}
	A holomorphic mapping $f\colon\mathbb{C}^p\rightarrow\mathbb{C}\mathbb{P}^n$ is said to be of {\sl maximal rank} if at some point $z\in\mathbb{C}^p$, the differential $df$ satisfies $\rank(df)=\min\{p,n\}$.	
\end{defi}

In the following statements, we set
\begin{equation}
\label{kappa(p,n) defi}
\kappa(p,n)=
\max\{n+1-p,1\}
=
\begin{cases}
n+1-p,&\text{if}\quad 1\leq p< n\\
1,&\text{if} \quad n\leq p.
\end{cases}
\end{equation}
Here are the statements of our main results.
\begin{namedthm*}{Theorem A}
	Let $p\geq 1,n\geq2$ be positive integers. Let $F$ be the Fermat hypersurface of degree $d$ in $\mathbb{C}\mathbb{P}^n$, defined by the homogeneous polynomial 
	$$
	Q(\omega)=\sum_{i=0}^n\omega_i^d,
	$$
	where $\omega=[\omega_0:\omega_1:\dots:\omega_n]$ is a homogeneous coordinate of $\mathbb{C}\mathbb{P}^n$.
	If $d>(n+1)\kappa(p,n-1)$, then the image of every holomorphic mapping $f\colon\mathbb{C}^p\rightarrow F$ of maximal rank lies in a linear subspace of dimension at most $\bigg[\dfrac{n-1}{2}\bigg]$.
\end{namedthm*}

\begin{namedthm*}{Theorem B}
	Let $p\geq 1,n\geq2$ be positive integers. Let $F$ be the Fermat hypersurface of degree $d$ in $\mathbb{C}\mathbb{P}^n$, defined by the homogeneous polynomial 
	$$
	Q(\omega)=\sum_{i=0}^n\omega_i^d,
	$$
	where $\omega=[\omega_0:\omega_1:\dots:\omega_n]$ is a homogeneous coordinate of $\mathbb{C}\mathbb{P}^n$.
	If $d>(n+1)\,\kappa(p,n)$, then the image of every holomorphic mapping $f\colon\mathbb{C}^p\rightarrow\mathbb{C}\mathbb{P}^n$ of maximal rank omitting the Fermat hypersurface $F$ lies in a linear subspace of dimension at most $\bigg[\dfrac{n}{2}\bigg]$.	
\end{namedthm*}

When $p\geq n$, the above two results are direct consequences of the Carlson-Griffiths' theory of equidimensional holomorphic mappings \cite{Carlson-Griffiths 72}. Thus Theorems A and B interpolate the works of Green and Carlson-Griffiths. Assume furthermore that $p>\bigg[\dfrac{n-1}{2}\bigg]$ in Theorem A or $p>\bigg[\dfrac{n}{2}\bigg]$ in Theorem B, we obtain the nonexistence of holomorphic mapping of maximal rank (Corollaries \ref{nonexistence compact case}, \ref{nonexistence logarithmic case}).

Etesse \cite{Etesse2023} proved that in the setting of Theorem A,  the image of the map $f$ must be contained in some Fermat hypersurface $F'$ in the subspace $\mathbb{C}\mathbb{P}^{n-1}$. One may think about applying successively this result to get more information about the degeneracy locus of $f$. However, this requires that the  induced map $\mathbb{C}^p\rightarrow F'$ must preserve the condition of maximal rank, which is  not always the case when  $p\geq 2$.

In \cite{Huynh2023}, using Nevanlinna theory, we show the degeneracy of the map $f$ in  the settings of the Theorems A and B. In the current paper, to overcome the technical problem of losing the maximal rank condition mentioned above,  we need to establish a Second Main Theorem for holomorphic mappings from $\mathbb{C}^p$ into $\mathbb{C}\mathbb{P}^n$ of {\em arbitrary rank}.

The paper is organized as follows. 
In Section~\ref{section: preparation}, we give a short introduction to higher dimensional Nevanlinna theory, and establish a Second Main Theorem type estimate for holomorphic mappings $f\colon\mathbb{C}^p\rightarrow\mathbb{C}\mathbb{P}^n$. Details of the proof of Theorem A and Theorem B will be provided in the section 3. In the last section, we provide an application of the two main results about the nonexistence of holomorphic mapping of maximal rank into or avoiding Fermat type hypersurfaces.
\section{Preliminaries}
\label{section: preparation}

\subsection{A brief introduction to higher dimensional Nevanlinna theory}
Let $z=(z_1,\dots,z_p)$ be the standard coordinate system of $\mathbb{C}^p$. Denote by $\|\cdot\|$ the Euclidean norm:
$$
\|z\|=\sqrt{\sum_{i=1}^{p}|z_i|^2}.
$$	
Set
$$
\alpha=\dif\!\dif^c\|z\|^2,\quad\beta=\dif\!\dif^c\log\|z\|^2,\quad\gamma=\dif^c\log \|z\|^2\wedge\beta^{p-1},
$$
where $\dif=\partial+\bar{\partial}$ and $\dif^c=\dfrac{i}{4\pi}(\bar{\partial}-\partial)$.

Let $\BB_r:=\{z\in\mathbb{C}^p:\|z\|< r\}\subset \mathbb{C}^p$ be the open ball in $\mathbb{C}^p$ of radius $r>0$ centered at the origin. Fix a truncation level $m\in \mathbb{N}\cup \{\infty\}$, for an effective divisor $E=\sum_i\alpha_i E_i$ on $\mathbb{C}^p$ where  $E_i$ are irreducible components and $\alpha_i\geq 0$,  the $m$-truncated degree of the divisor $E$ on  the balls is given by
\[
n^{[m]}(t,E)
:=
\dfrac{1}{t^{2p-2}}\int_{\mathbb{B}_t\cap(\sum_i\min\{m,\alpha_i\}E_i)}\alpha^{p-1}
\eqno
{{\scriptstyle (t\,>\,0)},}
\]
the \textsl{truncated counting function at level} $m$ of $E$ is then defined by taking the logarithmic average
\[
N^{[m]}(r,E)
\,
:=
\,
\int_1^r \frac{n^{[m]}(t, E)}{t}\,\dif\! t
\eqno
{{\scriptstyle (r\,>\,1)}.}
\]
When $m=\infty$, for abbreviation we  write $n(t,E)$, $N(r,E)$ for $n^{[\infty]}(t,E)$, $N^{[\infty]}(r,E)$ respectively.

Let $f\colon\mathbb{C}^p\rightarrow \mathbb{C}\mathbb{P}^n$ be an entire holomorphic curve having a reduced representation $f=[f_0:\cdots:f_n]$ in the homogeneous coordinates $[z_0:\cdots:z_n]$ of $\mathbb{C}\mathbb{P}^n$. Let $D=\{Q=0\}$ be a divisor in $\mathbb{C}\mathbb{P}^n$ defined by a homogeneous polynomial $Q\in\mathbb{C}[z_0,\dots,z_n]$ of degree $d\geq 1$. If $f(\mathbb{C}^p)\not\subset\supp D$, then $f^*D$ is a divisor on $\mathbb{C}^p$. We then define the \textsl{truncated counting function} of $f$ with respect to $D$ as
\[
N_f^{[m]}(r,D)
\,
:=
\,
N^{[m]}\big(r,f^*D\big),
\]
which measures the intersection frequency of $f(\mathbb{C}^p)$ with $D$. If $f^*D=\sum_i\mu_iE_i$, where $\mu_i>0
$ and $\mu=\min_i\{\mu_i\}$, then we say that  $f$ is {\sl completely $\mu$--ramified} over $D$, with the convention that $\mu=\infty$ if $f(\mathbb{C}^p)\cap\supp D=\varnothing$. Next,
the \textsl{proximity function} of $f$ associated to the divisor $D$ is given by
\[
m_f(r,D)
\,
:=
\,
\int_{\|z\|=r}
\log
\frac{\big\Vert f(z)\big\Vert_{\max}^d\,
	\Vert Q\Vert_{\max}}{\big|Q(f)(z)\big|}
\,
\gamma(z),
\]
where $\Vert Q\Vert_{\max}$ is the maximum  absolute value of the coefficients of $Q$ and where
\begin{equation}
\label{| |max definition}
\big\Vert f(z)\big\Vert_{\max}
:=
\max
\{|f_0(z)|,\dots,|f_n(z)|\}.
\end{equation}
Since $\big|Q(f)\big|\leq
\left(\substack{d+n\\ n}
\right)\,
\Vert Q\Vert_{\max}\cdot\Vert f\Vert_{\max}^d$, we see that $m_f(r,D)\geq O(1)$ is bounded  from below by some constant.
Lastly, the \textsl{Cartan order function} of $f$ is defined by
\begin{align*}
T_f(r)
:=
\int_{\|z\|=r}
\log
\big\Vert f(z)\big\Vert_{\max} \gamma(z).
\end{align*}

The Nevanlinna theory is then established by comparing the above three functions. It consists of two fundamental theorems (for  comprehensive expositions, see Noguchi-Winkelmann \cite{Noguchi-Winkelmann2014} and Ru \cite{Ru2021}).

\begin{namedthm*}{First Main Theorem}\label{fmt} Let $f\colon\mathbb{C}^p\rightarrow \mathbb{C}\mathbb{P}^n$ be a holomorphic map and let $D$ be a hypersurface of degree $d$ in $\mathbb{C}\mathbb{P}^n$ such that $f(\mathbb{C}^p)\not\subset D$. Then one has the estimate
	\[
	m_f(r,D)
	+
	N_f(r,D)
	\,
	=
	\,
	d\,T_f(r)
	+
	O(1)
	\]
	for every $r>1$,
	whence
	\begin{equation}
	\label{-fmt-inequality}
	N_f(r,D)
	\,
	\leq
	\,
	d\,T_f(r)+O(1).
	\end{equation}
\end{namedthm*}

Hence the First Main Theorem provides an upper bound on the counting function in term of the order function. The reverse direction, called {\sl Second Main Theorem}, is usually much harder, and one often needs to take the sum of the counting functions of many divisors.

Throughout this paper, for an entire holomorphic map $f$, the notation $S_f(r)$ means a real function of $r \in \mathbb{R}^+$ such that 
\[
S_f(r) \leq
O(\log(T_f(r)))+ \epsilon \log r
\]
for every positive constant $\epsilon$ and every $r$ outside of a subset (depending on $\epsilon$) of finite Lebesgue measure of $\mathbb{R}^+$. In the case where $f$ is rational, we understand that $S_f(r)=O(1)$. In any case we always have

\[
\liminf_{r\rightarrow\infty}\dfrac{S_f(r)}{T_f(r)}
=
0.
\]

\subsection{Generalized Wronskian}
Consider a collection of $n+1$ holomorphic functions $f_i\colon\mathbb{C}^p\rightarrow \mathbb{C}$ ($0\leq i\leq n$). For a tuple of $p$ nonnegative integers $\alpha=(\alpha_1,\dots,\alpha_p)$, let  $\Delta^{\alpha}$ denote the differential operator
\begin{equation}
\label{form of Delta alpha}
\Delta^{\alpha}
:=
\dfrac{\partial^{|\alpha|}}{\partial z_1^{\alpha_1}\dots\partial z_p^{\alpha_p}},
\end{equation}
having order $|\alpha|=\sum_{i=1}^{p}\alpha_i$. A collection $\mathcal{S}=\{\Delta^s\}_{0\leq s\leq n}$ of $n+1$ differential operators of the form \eqref{form of Delta alpha} is said to be {\sl admissible} if for each $0\leq s\leq n$, the order $|\Delta_s|:=i_1+\dots+i_p$ of
$$\Delta^s=
\dfrac{\partial^{i_1+\dots+i_p}}{\partial z_1^{i_1}\dots\partial z_p^{i_p}}$$   satisfies
$|\Delta_s|\leq s$. In particular when $s=0$, the operator $\Delta^0$ is the identity.

A generalized Wronskian of the family of holomorphic functions $\{f_0,\dots,f_n\}$ on $\mathbb{C}^p$ is the determinant of a square matrix of the form
$$
W_{\mathcal{S}}(f_0,f_1,\dots,f_n)
=
\det
\begin{pmatrix}
\Delta^0 (f_0)&\Delta^0{f_1}&\dots&\Delta^0(f_n)\\
\Delta^1(f_0)&\Delta^1(f_1)&\dots&\Delta^1(f_n)\\
\vdots&\vdots&\dots&\vdots\\
\Delta^n(f_0)&\Delta^n(f_1)&\dots&\Delta^n(f_n)
\end{pmatrix},
$$
where $\mathcal{S}=\{\Delta^s\}_{0\leq s\leq n}$ is an admissible family. When $p=1$, this boils down to the classical Wronskian. Similarly as in the one-dimensional case, generalized Wronskians can be employed to test the linearly dependence of the family $f_i$ \cite{Rot55,Fujimoto1985,Schmidt1980,BD2010}.

Now let $f\colon\mathbb{C}^p\rightarrow\mathbb{C}\mathbb{P}^n$ be a holomorphic mapping. If $\mathsf{s}=\max_{z\in\mathbb{C}^p}\rank(df_{z})$, then $f$ is  said to be  of rank $\mathsf{s}$. For this class of holomorphic mappings, we have the following result due to Fujimoto \cite{Fujimoto1985}.

\begin{thm}
	\label{generalized wronskian for rank s holo mapping Cp to CPn}
	Let $f\colon\mathbb{C}^p\rightarrow\mathbb{C}\mathbb{P}^n$ be a linearly nondegenerate holomorphic mapping of rank $\mathsf{s}$ and let $[f_0:f_1:\dots:f_n]$ be a reduced representation of $f$. Then there exists an admissible family $\mathcal{S}=\{\Delta^{\ell}\}_{0\leq \ell\leq n}$ containing at least $\mathsf{s}$ differential operators of order $1$ such that $W_{\mathcal{S}}(f_0,f_1,\dots,f_n)\not\equiv 0$. 
\end{thm}

\begin{rem}
Recently, Etesse \cite{Etesse2023} observes that in order to test the linearly dependence of the family of holomorphic functions $\{f_0,\dots,f_n\}$ on $\mathbb{C}^p$, it suffices to work on  a smaller subfamily of generalized Wronskians, called {\sl geometric generalized Wronskians}. Each element in an admissible family $\mathcal{S}$ is determined by a word written in the lexicographic order with the alphabet $\{1,\dots,p\}$. Identifying the family $\mathcal{S}$ with the set of its corresponding words, we say that $\mathcal{S}$ is a {\sl full set} if and only if it satisfies the following property:

\begin{center}
	if a word $s$ belongs to $\mathcal{S}$, then so does every one of its subwords.
\end{center}

A geometric generalized Wronskian of the family of $n+1$ holomorphic functions $\{f_0,\dots,f_n\}$ on $\mathbb{C}^p$ is a generalized Wronskian $W_{\mathcal{S}}(f_0,f_1,\dots,f_n)$, where $\mathcal{S}$ is a full set. Etesse  \cite[Theorem 1.4.1]{Etesse2023} showed that the the family $\{f_0,f_1,\dots,f_n\}$ is linearly independent over $\mathbb{C}$ if and only if there exists some geometric generalized Wronskian $$W_{\mathcal{S}}(f_0,f_1,\dots,f_n)$$  such that $W_{\mathcal{S}}(f_0,f_1,\dots,f_n)\not\equiv 0$. Using this observation, one can get another proof for Theorem~\ref{generalized wronskian for rank s holo mapping Cp to CPn} (see  \cite[Theorem 2.5]{Huynh2023} for more details).	
\end{rem}

\subsection{A Second Main Theorem for holomorphic mappings of arbitrary rank}

Applying Theorem~\ref{generalized wronskian for rank s holo mapping Cp to CPn}, together with the classical Wronskian method due to Cartan \cite{Cartan1933}, one obtains
\begin{thm}
	\label{smt for rank s} Let $1\leq p\leq n$ be two positive integers.
	Let $f:\mathbb{C}^p\rightarrow\mathbb{C}\mathbb{P}^n$ be a linearly nondegenerate  holomorphic mapping of  rank $s$, and let $\{H_i\}_{1\leq i\leq q}$ be a family of $q\geq n+2$ hyperplanes in general position in $\mathbb{C}\mathbb{P}^n$. Then the following Second Main Theorem type estimates holds
	\begin{align}
	\label{smt statement}
	(
	q
	-
	n-1
	)\,
	T_{f}(r)\leq\sum_{i=1}^q N_{f}^{[n+1-s]}(r,H_i)
	+
	S_{f}(r).
	\end{align}
\end{thm}

The reader is refereed to \cite{Huynh2023} for details of the proof. Note that in the case where $f$ is of maximal rank, such above Second Main Theorem was established by Noguchi \cite{Noguchi1996}.

It is well-known that one can deduce from a Second Main Theorem type estimate  the following two standard inequalities, called the defect relation and ramification theorem.
\begin{namedthm*}{Defect relation}
	Let $f:\mathbb{C}^p\rightarrow\mathbb{C}\mathbb{P}^n$ be a linearly nondegenerate  holomorphic mapping of rank $\mathsf{s}$, and let $\{H_i\}_{1\leq i\leq q}$ be a family of $q\geq n+2$ hyperplanes in general position in $\mathbb{C}\mathbb{P}^n$. Then one has the following defect relation
	\begin{equation}
	\label{defect relation}
	\sum_{i=1}^q \delta_{f}^{[\kappa(s,n)]}(H_i)
	\leq
	n+1,
	\end{equation}
	where $\kappa(s,n)$ is given as in~\eqref{kappa(p,n) defi}
\end{namedthm*}
For a divisor $D$ on $\mathbb{C}\mathbb{P}^n$, if $f^*D=\sum_j\mu_jA_j$ is the decomposition into irreducible components and $\mu=\min_j\mu_j$, then $f$ is said to be {\em completely $\mu$--ramified} over $D$, with the convention that $\mu=\infty$ when $f^{-1}D=\varnothing$.
\begin{namedthm*}{Ramification Theorem}
	Let $f:\mathbb{C}^p\rightarrow\mathbb{C}\mathbb{P}^n$ be a linearly nondegenerate holomorphic mapping of rank $\mathsf{s}$, and let $\{H_i\}_{1\leq i\leq q}$ be a family of $q\geq n+2$ hyperplanes in general position in $\mathbb{C}\mathbb{P}^n$. If $f$ is completely $\mu_{i}$--ramified over each  hyperplane $H_i$ for $i=1,\dots, q$, then one has
	\[
	\label{ramification theorem statement}
	\sum_{i=1}^q
	\bigg(
	1
	-
	\dfrac{\kappa(s,n)}{\mu_{i}}
	\bigg) 
	\leq
	n+1,
	\]
	where $\kappa(s,n)$ is given as in \eqref{kappa(p,n) defi}.
\end{namedthm*}

\section{Proof of the main results}
\subsection{Proof of Theorem A}
\begin{proof}
	We follow the arguments of \cite[Example 3.10.21]{Kob98} (see also \cite{Green75}). Let $H\cong\mathbb{C}\mathbb{P}^{n-1}$ be the hyperplane in $\mathbb{C}\mathbb{P}^n$ defined by $\sum_{i=0}^{n}\omega_i=0$. Let $[f_0:f_1:\dots:f_n]$ be a reduced representation of $f$. Consider the endomorphism $$\pi\colon\mathbb{C}\mathbb{P}^n\rightarrow\mathbb{C}\mathbb{P}^n,\qquad [\omega_0:\omega_1:\dots:\omega_n]\rightarrow [\omega_0^d:\omega_1^d:\dots:\omega_n^d],
$$
and put $g=\pi\circ f\colon\mathbb{C}^p\rightarrow H\cong\mathbb{C}\mathbb{P}^{n-1}$. Then $g$ is of maximal rank and its image lies in the hyperplane $H\cong\mathbb{C}\mathbb{P}^{n-1}$. Let $\{H_i\}_{0\leq i\leq n}$ be the family of $n+1$ hyperplanes in $H\cong\mathbb{C}\mathbb{P}^{n-1}$ given by $H_i=\{\omega_i=0\}$, which is in general position. For any $z\in g^{-1}(H_i)=f^{-1}(H_i)$, it is clear that $\ord_zg^*H_i\geq d$. If the image of $g$ doesn't lie in a smaller linear subspace of $H$, then by applying the ramification theorem for $g$ and the family $\{H_i\}_{0\leq i\leq n}$, one obtains
\[
\sum_{i=0}^{n}\bigg(1-\dfrac{\kappa(p,n-1)}{d}\bigg)\leq n,
\]
which yields $d\leq (n+1)\kappa(p,n)$, a contradiction. Hence the image of $g$ must be contained in some smaller linear subspace of $H$, which implies that $f_i$ must satisfy another nontrivial  relation
\begin{equation}
\label{second linear equation f_i satisfies}
\sum_{i=0}^{n}a_if_i^d=0,
\end{equation}
where $(a_0,a_1,\dots,a_n)\not\in\{(1,1,\dots,1),(0,0,\dots.0)\}$. Without lost of generality, one may assume that $a_n=1$. Then one has
\begin{equation}
\label{equation from the first and second for fi}
\sum_{i=0}^{n-1}(a_i-1)f_i^d=0.
\end{equation}
Suppose that $f_i\not\equiv0$ for any $1\leq i\leq n-1$, otherwise the problem reduces to lower dimensional case. We first prove the following
\begin{claim}
	There exist some indexes $i\not=j$ with $1\leq i,j\leq n-1$ such that $f_i/f_j$ is constant.
\end{claim}
Indeed, put $I=\{i:0\leq i\leq n-1,a_i\not=1\}$, then $2\leq |I|\leq n$ and \eqref{equation from the first and second for fi} can be rewritten as
$\sum_{i\in I}c_if_i=0$, where $c_i\not=0$. If $|I|=2$, we finish the proof. Suppose $|I|=\gamma\geq3$. Consider the homogeneous coordinates $[\omega_i]_{i\in I}$ of the linear subspace $\mathbb{C}\mathbb{P}^{\gamma-1}$. Similar as in above, consider the hyperplane $H_I\cong \mathbb{C}\mathbb{P}^{\gamma-2}\subset \mathbb{C}\mathbb{P}^{\gamma-1}$ defined as $\sum_{i\in I}c_i\omega_i=0$ and the holomorphic mapping 
$$
\mathsf{f}_{I}\colon\mathbb{C}^p\rightarrow \mathbb{C}\mathbb{P}^{\gamma-1},\qquad \mathsf{f}(z)=[f_i(z)]_{i\in I}.
$$ Denote by $\pi_{I}$ the endomorphism $\pi_{I}\colon \mathbb{C}\mathbb{P}^{\gamma-1}\rightarrow \mathbb{C}\mathbb{P}^{\gamma-1},\,[\omega_i]\rightarrow[\omega_i^d]$. Then the image of $g_I:=\pi_{I}\circ \mathsf{f}_I$ lies in  $H_I$. Put $s=\rank(g_I)$, then
\[
s\geq\max\{ p-(n+1-\gamma),1\}.
\] If the image of $g_I$ doesn't lie in some smaller linear subspace, then using ramification theorem for $g_I$ and the family of $\gamma$ hyperplanes $\{H_i\}_{i\in I}$ in $H_I\cong \mathbb{C}\mathbb{P}^{\gamma-2}$, one obtains 
\begin{equation}
\label{applying ramification for gI}
\sum_{i\in I}\bigg(1-\dfrac{d_I}{d}\bigg)\leq \gamma-1,
\end{equation}
where $d_I=\max\{\gamma-s-1,1\}$. This implies
$d\leq \gamma d_I$. If $p+\gamma>n+1$, then $s=p+\gamma-n-1$, hence $\gamma-s-1=n-p$. In this case $\gamma d_I=\gamma(n-p)<(n+1)(n-p)$, a contradiction. In the case where $p+\gamma\leq n+1$, one has $s=1$ and $d_I=\gamma-2$. Hence $\gamma d_I=\gamma(\gamma-2)\leq (n+1-p)(n-p-1)<(n+1)(n-p)<d$, contradiction. Thus $g_I$ is linearly degenerate. Inductively, the claim is proved.

Going back to the proof of the theorem, we now follow the arguments in \cite[Example 3.10.21]{Kob98}. Let $\sim$ be the equivalence relation on the index set $\{0,1,\dots,n\}$ defined as $i\sim j$ if and only if $f_i/f_j$ is constant and let $\{I_1,\dots, I_m\}$ be the partition of $\{0,1,\dots,n\}$ by $\sim$. Then for each $1\leq r\leq m$, we pick an index $i_r\in I_r$ and put
$f_j=\ell_jf_{i_r}$. Set $b_r=\sum_{j\in I_r}\ell_j$, then
\[
\sum_{r=1}^mb_rf_{i_r}^d=0.
\]
Using the above claim, one gets $b_r=0$ for all $1\leq r\leq m$. Hence each $I_r$ contains at least $2$ indexes. Thus the image of $f$ lies in the linear subspace defined by the equations
\begin{equation}
\label{linear relation f_i}
\omega_j=\ell_j\omega_{i_r},\quad 1\leq r\leq m,\quad j\in I_r,\, j\not=i_r.
\end{equation}
Consequently, $f_i$ must satisfy at least $\sum_{r=1}^m( |I_r|-1)=n+1-m$ independent linear relations. Since $|I_r|\geq 2 $ for $1\leq r\leq m$, it is trivial that $m\leq \bigg[\dfrac{n+1}{2}\bigg]$, which implies that the image of $f$ is contained in a linear subspace of dimension at most $\bigg[\dfrac{n+1}{2}\bigg]-1=\bigg[\dfrac{n-1}{2}\bigg]$.
\end{proof}

\subsection{Proof of Theorem B}

\begin{proof}
 Suppose that $[f_0:f_1:\dots:f_n]$ is a reduced representation of $f$. Since $f$ avoids the Fermat hypersurface $F$, one can write
 $$f_0^d+f_1^d+\dots+f_n^d
 =e^h,$$
 for some entire holomorphic map $h$. Set $f_{n+1}=e^{\frac{h+i\pi}{d}}$, then the above equality can be rewritten as $$f_0^d+f_1^d+\dots+f_{n+1}^d=0.
 $$
Consider the holomorphic map $\mathsf{f}\colon\mathbb{C}^p\rightarrow\mathbb{C}\mathbb{P}^{n+1}$, given by $[f_0:f_1:\dots:f_{n+1}]$. Then $\mathsf{f}$ is also of maximal rank and its image is contained in the Fermat hypersurface of degree $d$ in $\mathbb{C}\mathbb{P}^{n+1}$. As in the proof of Theorem A, let $H\cong\mathbb{P}^{n}(\mathbb{C})$ be the hyperplane in $\mathbb{P}^{n+1}(\mathbb{C})$ defined by $\sum_{i=0}^{n+1}\omega_i=0$, consider the endomorphism $$\pi\colon\mathbb{P}^{n+1}(\mathbb{C})\rightarrow\mathbb{P}^{n+1}(\mathbb{C}),\qquad [\omega_0:\omega_1:\dots:\omega_{n+1}]\rightarrow [\omega_0^d:\omega_1^d:\dots:\omega_{n+1}^d],
$$
and put $\mathsf{g}=\pi\circ \mathsf{f}\colon\mathbb{C}^p\rightarrow\mathbb{P}^{n+1}(\mathbb{C})$. Then $\mathsf{g}$ is also of maximal rank. Let $\{H_i\}_{0\leq i\leq n+1}$ be the family of $n+2$ hyperplanes in $H\cong\mathbb{C}\mathbb{P}^{n+1}$ given by
\begin{align*}
H_i&=\{\omega_i=0\},\quad(0\leq i\leq n+1),
\end{align*}
 which is in general position. If the image of $\mathsf{g}$ doesn't lie in any hyperplane, then by applying the ramification theorem for $\mathsf{g}$ and the family $\{H_i\}_{0\leq i\leq n+1}$, noting that $\mathsf{g}$ omits $H_{n+1}$, one obtains
	\[
	\sum_{i=0}^{n}\bigg(1-\dfrac{\kappa(p,n)}{d}\bigg)+1\leq n+1,
	\]
	which yields $d\leq (n+1)\kappa(p,n)$, a contradiction. Hence the image of $\mathsf{g}$ must be contained in some hyperplane. Continuing the arguments of the proof of Theorem A, we can find a partition $\{I_1,\dots, I_m\}$  of $\{0,1,\dots,n+1\}$ and indexes $i_r\in I_r$ such that $|I_r|\geq 2$ and the image of $\mathsf{f}$ lies in the linear subspace defined by the equations
	$$\omega_j=\ell_j\omega_{i_r},\quad 1\leq r\leq m,\quad j\in I_r,\, j\not=i_r.$$
If $n+1\in I_{r_0}$, then by changing the index $i_{r_0}$ if necessary, one may assume $i_{r_0}\not=n+1$ and hence, after removing the equation $\omega_{n+1}=\ell_{n+1}\omega_{i_{r_0}}$, the family $\{f_i\}_{0\leq f\leq n}$ must satisfy at least $\sum_{r=1}^m( |I_r|-1)-1=n+1-m$ independent linear relations. Since $|I_r|\geq 2$ for $1\leq r\leq m$, there holds $m\leq \bigg[\dfrac{n+2}{2}\bigg]$, which implies that the image of $f$ is contained in a linear subspace of dimension at most $\bigg[\dfrac{n+2}{2}\bigg]-1=\bigg[\dfrac{n}{2}\bigg]$.
\end{proof}	

\section{Some discussions}
Keeping the notations as in the proof of Theorem A, then  one can find at least $n+1-m$ independent linear relations among the column vectors of  the  Jacobian of $f$. Consequently, one has
\[
\rank(df)\leq 
\bigg[\dfrac{n-1}{2}\bigg].
\]
Therefore if $p>\bigg[\dfrac{n-1}{2}\bigg]$, there exists no holomorphic mapping of maximal rank $f\colon\mathbb{C}^p\rightarrow F\subset\mathbb{C}\mathbb{P}^n$. Similar result holds true in the setting of Theorem B.
\begin{cor}
	\label{nonexistence compact case}
	Let $p\geq 1,n\geq2$ be positive integers. Let $F$ be the Fermat hypersurface of degree $d$ in $\mathbb{C}\mathbb{P}^n$, defined by the homogeneous polynomial 
	$$
	Q(\omega)=\sum_{i=0}^n\omega_i^d,
	$$
	where $\omega=[\omega_0:\omega_1:\dots:\omega_n]$ is a homogeneous coordinate of $\mathbb{C}\mathbb{P}^n$.
	If $d>(n+1)\kappa(p,n-1)$ and $p>\bigg[\dfrac{n-1}{2}\bigg]$, then there exists no holomorphic mapping  $f\colon\mathbb{C}^p\rightarrow F$ of maximal rank.
\end{cor}

\begin{cor}
\label{nonexistence logarithmic case}
	Let $p\geq 1,n\geq2$ be positive integers. Let $F$ be the Fermat hypersurface of degree $d$ in $\mathbb{C}\mathbb{P}^n$, defined by the homogeneous polynomial 
	$$
	Q(\omega)=\sum_{i=0}^n\omega_i^d,
	$$
	where $\omega=[\omega_0:\omega_1:\dots:\omega_n]$ is a homogeneous coordinate of $\mathbb{C}\mathbb{P}^n$.
	If $d>(n+1)\kappa(p,n)$ and $p>\bigg[\dfrac{n}{2}\bigg]$, then there exists no holomorphic mapping $f\colon\mathbb{C}^p\rightarrow\mathbb{C}\mathbb{P}^n\setminus F$ of maximal rank.
\end{cor}
Note that in \cite{Pacienza-Rousseau2012}, by using the jet-bundle technique, Pacienza-Rousseau proved the non-existence of  maximal rank holomorphic maps from $\mathbb{C}^2$ into the very general degree $d$ hypersurface of $\mathbb{C}\mathbb{P}^4$ for high enough degree $d\geq 93$. This result may be seen as first steps towards the $p$-measure hyperbolicity of generic projective hypersurfaces of high degree and provides support evident to the following
\begin{namedthm*}{Generalized Kobayashi's conjecture }
Let $X_d$ be a generic hypersurface of degree $d$ in $\mathbb{CP}^n$.
\begin{itemize}
	\item If $d\geq 2n+1-p$ then $X_d$ is $p$-measure hyperbolicity ($1\leq p\leq n-1$).
	\item If $d\geq 2n+2-p$ then $X_d$ is $p$-measure hyperbolicity ($1\leq p\leq n$).
\end{itemize}	
\end{namedthm*}

Corollaries~\ref{nonexistence compact case}, \ref{nonexistence logarithmic case} provide further evidents to the above conjecture. Finally, we close this paper by remaining the following
\begin{namedthm*}{Conjecture}
The truncation level in Theorem \ref{smt for rank s} can be decreased to $1$, provided that  $f\colon\mathbb{C}^p\rightarrow\mathbb{C}\mathbb{P}^n$ is algebraically nondegenerate.
\end{namedthm*}
This expected statement would yield the Green-Griffiths conjecture for the case of Fermat type hypersurfaces, see \cite{HV2021} for supported evidents.

\section*{Funding}
This research is funded by University of Education, Hue University under grant number NCTB-T.24-TN.101.01.

\begin{center}
	\bibliographystyle{plain}
	
\end{center}
\address
\end{document}